\documentclass[journal,twoside,web]{ieeecolor}
\usepackage{lcsys}
\usepackage{cite}
\usepackage{amsmath,amssymb,amsfonts}

\usepackage{graphicx}
\usepackage{textcomp}

\usepackage{mathtools}
\usepackage{amsmath} 
\usepackage{amssymb}  

\usepackage{units}
\let\labelindent\relax
\usepackage{enumitem}
\usepackage{hyperref}
\usepackage{float}
\usepackage{algorithm}
\usepackage{algpseudocode}
\usepackage{cite}
\newtheorem{definition}{Definition}
\newtheorem{remark}{Remark}
\newtheorem{theorem}{Theorem}
\newtheorem{lemma}{Lemma}
\newtheorem{problem}{Problem}

\newtheorem{assumption}{Assumption}

\newcommand{\defeq}{\vcentcolon=}
\newcommand{\T}{\scriptscriptstyle\top}       
\newcommand{\mathst}{\text{s.t.}}

\DeclareMathOperator*{\argmax}{arg\,\operatorname*{max}}

\DeclareUnicodeCharacter{2212}{-}

\allowdisplaybreaks

\makeatletter
\def\widebreve{\mathpalette\wide@breve}
\def\wide@breve#1#2{\sbox\z@{$#1#2$}%
     \mathop{\vbox{\m@th\ialign{##\crcr
\kern0.18em\vspace{-0.016cm}\brevefill#1{0.55\wd\z@}\crcr\noalign{\nointerlineskip}%
                    $\hss#1#2\hss$\crcr}}}\limits}
\def\brevefill#1#2{$\m@th\sbox\tw@{$#1($}%
  \hss\resizebox{#2}{\wd\tw@}{\rotatebox[origin=c]{90}{\upshape(}}\hss$}
\makeatletter

\def\BibTeX{{\rm B\kern-.05em{\sc i\kern-.025em b}\kern-.08em
    T\kern-.1667em\lower.7ex\hbox{E}\kern-.125emX}}
\usepackage{eso-pic}
\AddToShipoutPictureBG*{%
\AtPageUpperLeft{%
\setlength\unitlength{1in}%
\hspace*{\dimexpr0.5\paperwidth\relax}
\makebox(0,-0.55)[c]{
\begin{tabular}{c c}

To appear in the IEEE Control Systems Letters.

\end{tabular}}}}

\begin{document}
\title{Counterexample-Guided Synthesis of Robust Discrete-Time Control Barrier Functions}
\author{Erfan Shakhesi, \IEEEmembership{Student Member, IEEE}, Alexander Katriniok, \IEEEmembership{Senior Member, IEEE}, and \mbox{W.P.M.H. (Maurice) Heemels}, \IEEEmembership{Fellow, IEEE}
\thanks{Received March 17, 2025; revised May 18, 2025; accepted \mbox{June 2,
2025}. Date of publication XXX XX 2025, date of current version XXX
XX 2025. The research has received funding from the European Research Council under the Advanced ERC Grant Agreement PROACTHIS, \mbox{no. 101055384.}}
\thanks{The authors are with the Control Systems Technology \mbox{section}, Mechanical Engineering, Eindhoven University of Technology, The Netherlands. E-mail: {\{\tt\small e.shakhesi, a.katriniok, m.heemels\}@tue.nl}
}}

\maketitle

\begin{abstract}
Learning-based methods have gained popularity for training candidate Control Barrier Functions (CBFs) to satisfy the CBF conditions on a finite set of sampled states. However, since the CBF is unknown \textit{a priori}, it is unclear which sampled states belong to its zero-superlevel set and must satisfy the CBF conditions, and which ones lie outside it. Existing approaches define a set in which all sampled states are required to satisfy the CBF conditions, thus introducing conservatism. In this paper, we address this issue for robust discrete-time CBFs (R-DTCBFs). Furthermore, we propose a class of R-DTCBFs that can be used in an online optimization problem to synthesize safe controllers for general discrete-time systems with input constraints and bounded disturbances. To train such an R-DTCBF that is valid not only on sampled states but also across the entire region, we employ a verification algorithm iteratively in a counterexample-guided approach. We apply the proposed method to numerical case studies.
\end{abstract}
\begin{IEEEkeywords}
Control Barrier Functions, Lyapunov methods, neural networks, robust control.
\end{IEEEkeywords}
\section{Introduction}
\label{sec:introduction}
\IEEEPARstart{C}{ontrol} Barrier Functions (CBFs) are Lyapunov-like functions that guarantee the controlled invariance of their zero-superlevel set, defined as the set of all states for which a CBF is non-negative \cite{Ames2014a}. Such a function can be used in an online optimization problem to synthesize safe controllers by adjusting known nominal controllers. CBFs are applied both to continuous-time systems \cite{Ames2014a} and discrete-time systems \cite{DCBF-Definition}, where, in the latter case, they are called discrete-time CBFs (DTCBFs). DTCBFs offer potential for various applications, including their use within Model Predictive Control frameworks \cite{katriniok2023}.
Moreover, to guarantee safety under disturbances, the definitions of CBFs and DTCBFs have been extended to robust CBFs \cite{JANKOVIC2018359} and robust DTCBFs (R-DTCBFs) \cite{RDTCBF-literature}. However, a key limitation of existing R-DTCBF formulations in the literature is that the associated constraint cannot be used in online optimization problems for safe controller synthesis due to the issue of an infinite number of constraints.

To synthesize robust (DT)CBFs, learning-based methods provide an attractive approach for training candidate robust (DT)CBFs based on a finite set of sampled states, as these methods are applied successfully in other settings, including training candidate Control Lyapunov Functions (CLFs) \cite{Neural-Lyapunov} and candidate Control Barrier Certificates (CBCs) \cite{Formally-Verified-NN}. Note that the term ``candidate'' indicates that the trained function is not necessarily valid over the entire region, but only on sampled states.
An additional and crucial difficulty in training candidate (DT)CBFs is that the CBF conditions are imposed only on the states within the zero-superlevel set of the resulting (DT)CBF, but this set is unknown \textit{a priori}. 
To bypass this issue, \cite{CBF-NN} considers an explicit region in which \textit{all} sampled states are required to satisfy the CBF conditions. This introduces \mbox{conservatism}, and in some applications, finding such a region may not be possible. Besides this work, \cite{labeling} employs an iterative approach to approximate the maximal region in which the CBF conditions can be satisfied. However, not all states within this approximated region satisfy the conditions, and thus the trained function is not necessarily a valid (DT)CBF. 



To ensure that a trained function is not only valid on \mbox{sampled} states but also on the entire region, two approaches are used in the literature: (i) sufficient sampling of the state space, such as \cite{Formally-Verified-NN} for CBCs, and \mbox{(ii) counterexample-guided} approaches, such as \cite{Neural-Lyapunov} for CLFs and \cite{CBF-NN} for CBFs. 
\mbox{Approach (i)} determines the maximum permissible distance between sampled states using Lipschitz constants and \mbox{applies} a threshold to ensure that satisfying the conditions at these states implies satisfaction throughout the region between them. \mbox{Although} conceptually elegant, it is computationally demanding with many samples and introduces conservatism. In approach (ii), a verification algorithm searches for a state that does not satisfy the required conditions across the entire region, not just at sampled states. This state is referred to as a counterexample. If a counterexample is found, it is added to the set of sampled states, and the candidate function is retrained. This procedure repeats until no counterexample is found. To the best of our knowledge, a counterexample-guided approach has not been applied to \mbox{(R-)DTCBFs}, likely due to the challenge of verifying candidate \mbox{(R-)DTCBFs} with trained neural network controllers in the discrete-time setting, where the (R-)DTCBF constraint is non-affine in the control input.

To overcome the above-mentioned challenges, our main contributions in this paper are: 
\begin{itemize}[leftmargin= 11pt]
    \item We introduce a class of R-DTCBFs that can be used in an online optimization problem to synthesize safe controllers for systems with bounded disturbances (Section \ref{sec:definition-RDCBF}).
    \item We extend the verification algorithm for DTCBFs \mbox{introduced} in our previous work \cite{Shakhesi2024} by modifying its steps so that it either verifies that a candidate R-DTCBF \mbox{belongs} to the proposed class or falsifies it by providing a \mbox{counterexample} within a predefined tolerance (Section~\ref{sec:RDCBF-Verification}). Our verification algorithm does not require a control policy, making it well-suited for a counterexample-guided approach to learning valid (R-)DTCBFs with less conservatism.
    \item We design a loss function that reduces conservatism in \mbox{training} candidate (DT)CBFs compared to existing \mbox{methods}, promotes larger zero-superlevel sets, and enforces the R-DTCBF conditions on \textit{designated} sampled states. We then propose a counterexample-guided approach to synthesize valid R-DTCBFs from the proposed class (Section \ref{sec:RDCBF-Train}).

    \item We apply our proposed methods to various numerical case studies (Section \ref{sec:numerical-case-studies}). 
\end{itemize}
\textit{Notation:} 
For a vector \mbox{$x \in \mathbb{R}^n$}, $x_i \in \mathbb{R}$, $i \in \{1,2, \hdots, n\}$, represents the \mbox{$i$-th} element of $x$, and $\lVert x\rVert$ is its Euclidean norm. 
A continuous function \mbox{$\gamma:\mathbb{R}_{\geqslant 0} \rightarrow \mathbb{R}_{\geqslant 0}$} is said to be a \textit{class-$\mathcal{K}_\infty$ function}, denoted by $\gamma \in \mathcal{K}_{\infty}$, if it is strictly increasing, $\gamma(0) = 0$, and $\lim_{r\rightarrow \infty} \gamma(r)=\infty$. For $\gamma \in \mathcal{K}_{\infty}$, we utilize $\gamma \in \mathcal{K}^{\leqslant \mathrm{id}}_{\infty}$ to indicate that $\gamma(r) \leqslant r$ for all $r \in \mathbb{R}_{\geqslant 0}$. A function $h:\mathbb{R}^n \rightarrow \mathbb{R}$ is \textit{Lipschitz continuous} on a set $\mathcal{C} \subseteq \mathbb{R}^n$, if there exists a \textit{Lipschitz constant} $L_h \in \mathbb{R}_{\geqslant 0}$ such that for all $x,y \in \mathcal{C}$, $|h(x) - h(y)| \leqslant L_h\lVert x - y\rVert$.
\section{Robust Discrete-time Control Barrier Function} \label{sec:definition-RDCBF}
We consider discrete-time systems of the form \vspace{-0.01cm}
\begin{align} \label{eq:discrete-time-dynamical-system}
    x_{t+1} = f(x_t, u_t) + w_t,
\end{align} 
where $x_t \in \mathbb{R}^n$ is the state vector, $u_t \in \mathbb{U} \subseteq \mathbb{R}^m$ is the control input vector, and $w_t \in \mathbb{W} \subset \mathbb{R}^n$ is the disturbance vector, all at the time instant $t \in \mathbb{N}_0$, and \mbox{$f:\mathbb{R}^n\times \mathbb{R}^m \rightarrow \mathbb{R}^n$} describes the unperturbed dynamics. Here, $\mathbb{U}$ is the control admissible set, and $\mathbb{W}$ is the disturbance set, which is bounded, i.e., there is $w_{\text{max}} \in \mathbb{R}_{\geqslant 0}$ such that for all $w \in \mathbb{W}$, $\lVert w \rVert \leqslant w_{\text{max}}$.
For the system \eqref{eq:discrete-time-dynamical-system}, the safe set $\mathcal{S}$ is defined as \vspace{-0.02cm}
\begin{align} \label{eq:safe-set}
    \mathcal{S} \defeq \{ x \in \mathbb{R}^{n} \mid s_i(x) \geqslant 0, ~ i \in \{1,2,\hdots, n_s\} \},
\end{align} \\ [-0.45cm]
where $s_i:\mathbb{R}^{n} \rightarrow \mathbb{R}$, $i \in \{1,2,\hdots, n_s\}$, are given mappings. 
\begin{definition}[Robust Controlled Invariance \cite{BLANCHINI1999}] \mbox{For the} system \eqref{eq:discrete-time-dynamical-system} with the control admissible set $\mathbb{U}$ and the disturbance set $\mathbb{W}$, a set $\mathcal{C} \subset \mathbb{R}^n$ is \textit{robust controlled invariant}, if, for each $x \in \mathcal{C}$, there exists a control input $u \in \mathbb{U}$ such that for all $w \in \mathbb{W}$, it holds that $f(x,u)+w \in \mathcal{C}$.
\end{definition}
\begin{definition}[Robust Safety] \label{def:safety}
    The system \eqref{eq:discrete-time-dynamical-system} with $\mathbb{U}$ and $\mathbb{W}$ is considered \textit{robustly safe} with respect to a safe set \mbox{$\mathcal{S} \subseteq \mathbb{R}^n$} and an initial state set $X_0 \subseteq \mathcal{S}$, denoted by $(\mathcal{S}, X_0)$-\textit{safe}, if there exists a control policy $\mu:\mathcal{S} \rightarrow \mathbb{U}$ such that for all trajectories of $x_{t+1} = f(x_t, \mu(x_t))+w_t$ with $x_0 \in X_0$ and $w_t \in \mathbb{W}$, it holds that $x_t \in \mathcal{S}$ for all $t \in \mathbb{N}_0$.
\end{definition}

To guarantee robust safety for the system \eqref{eq:discrete-time-dynamical-system}, one can construct a robust controlled invariant set, preferably sufficiently large, within the safe set $\mathcal{S}$.
To do so, the definition of discrete-time Control Barrier Functions \cite{DCBF-Definition} can be extended to robust discrete-time Control Barrier Functions, as in \cite{RDTCBF-literature}. 
\begin{definition}[R-DTCBF] \label{def:DTCBF}
    \hspace{-0.12cm}Consider a function $\hspace{-0.05cm}h:\hspace{-0.1cm} \mathbb{R}^n \hspace{-0.1cm}\rightarrow \hspace{-0.1cm} \mathbb{R}$ with zero-superlevel set \mbox{$\mathcal{C} \defeq \{x \in \mathbb{R}^n \mid h(x) \geqslant 0\}$}. For the system \eqref{eq:discrete-time-dynamical-system} with the control admissible set $\mathbb{U}$ and the disturbance set $\mathbb{W}$, $h$ is a \textit{robust discrete-time Control Barrier Function (R-DTCBF)}, if there exists a $\gamma \in \mathcal{K}^{\leqslant \mathrm{id}}_{\infty}$ such that for each $x \in \mathcal{C}$, there exists a control input $u \in \mathbb{U}$ that satisfies \vspace{-0.055cm}
    \begin{align} \label{eq:DTCBF-constraint-initial}
         h(f(x,u) + w)  - h(x) + \gamma(h(x)) \geqslant 0 \text{\hspace{0.18cm}for all  } w \in \mathbb{W}.
    \end{align} \\[-0.5cm]
    We also say that a control policy $\pi:\mathbb{R}^n \rightarrow \mathbb{R}^m$ is a \textit{friend} of an \textit{R-DTCBF pair} $(h, \gamma)$, if for each $x \in \mathcal{C}$, $\pi(x) \in \mathbb{U}$ and $h(f(x,\pi(x)) + w) - h(x) + \gamma(h(x)) \geqslant 0$ for all $w \in \mathbb{W}$.
\end{definition}

    Given that the disturbance $w_t$ at $t \in \mathbb{N}_0$ is unknown, the constraint \eqref{eq:DTCBF-constraint-initial} cannot be directly used in an online optimization problem to synthesize a safe controller due to the issue of an infinite number of constraints. To address this, we propose a class of R-DTCBFs in Lemma \ref{lemma:sufficient-RDTCBF}, which provides a sufficient condition for the satisfaction of \eqref{eq:DTCBF-constraint-initial}.
    
\begin{lemma} \label{lemma:sufficient-RDTCBF}
Let $h:\mathbb{R}^n \rightarrow \mathbb{R}$ be a Lipschitz continuous function on its zero-superlevel set $\mathcal{C}$. Suppose there exists a $\gamma \in \mathcal{K}_{\infty}^{\leqslant \text{id}}$ such that, for each $x \in \mathcal{C}$, there exists a control input $u \in \mathbb{U}$ satisfying the R-DTCBF constraint \vspace{-0.055cm}
\begin{align} \label{eq:DTCBF-constraint}
    \mathcal{H} (x,u) &\defeq h(f(x,u)) - h(x) + \gamma(h(x)) \geqslant L_hw_{\text{max}},
\end{align} \\[-0.5cm]
where $L_h \in \mathbb{R}_{\geqslant 0}$ is a Lipschitz constant for $h$ on $\mathcal{C}$. Then, $h$ is an R-DTCBF, and $\mathcal{C}$ is robust controlled invariant. Moreover, the system \eqref{eq:discrete-time-dynamical-system} is $(\mathcal{S}, \mathcal{C})$-safe, if $\mathcal{C} \subseteq \mathcal{S}$. \\
\begin{proof}
By the definition of Lipschitz continuity, for all $x \in \mathcal{C}$, $u \in \mathbb{U}$, and for all $w \in \mathbb{W}$, it follows that  \vspace{-0.075cm}
    \begin{align} \label{eq:DTCBF-proof-1}
        L_hw_{\text{max}} \geqslant L_h\lVert&w\rVert \geqslant h(f(x,u)) - h(f(x,u) + w). 
    \end{align} \\ [-0.51cm]
    By virtue of \eqref{eq:DTCBF-constraint} and \eqref{eq:DTCBF-proof-1}, there exists a $\gamma \in \mathcal{K}^{\leqslant \text{id}}_{\infty}$ such that for each $x \in \mathcal{C}$, an input $u \in \mathbb{U}$ exists that for all $w \in \mathbb{W}$ satisfies \vspace{-0.54cm}
    \begin{align} 
        h(f(x,u) + w) \geqslant h(x) - \gamma(h(x)) \geqslant 0. 
    \end{align} \\ [-0.52cm]
    Hence, $h$ is an R-DTCBF according to Definition \ref{def:DTCBF}. Moreover, this implies that if $x \in \mathcal{C}$, there exists an input $u \in \mathbb{U}$ such that $f(x, u)+w \in \mathcal{C}$ for all $w \in \mathbb{W}$, and thus $\mathcal{C}$ is \textit{robust controlled invariant}. Lastly, following \mbox{Definition \ref{def:safety}}, if $\mathcal{C} \subseteq \mathcal{S}$, the system \eqref{eq:discrete-time-dynamical-system} is $(\mathcal{S}, \mathcal{C})$-safe.
\end{proof}
\end{lemma}

Given an R-DTCBF pair $(h, \gamma)$ that satisfies Lemma \ref{lemma:sufficient-RDTCBF} for system \eqref{eq:discrete-time-dynamical-system} with $\mathbb{U}$ and $\mathbb{W}$, a nominal controller $\pi_{\text{nom}} \hspace{-0.06cm}:\hspace{-0.06cm}\mathbb{R}^n \hspace{-0.06cm}\rightarrow \hspace{-0.06cm}\mathbb{R}^m \hspace{-0.1cm}$, and $x_t$ at $t \in \mathbb{N}_0$, a safe control \mbox{input $u^*_t$} at this time instant can be obtained by solving the online optimization problem \vspace{-0.15cm}
\begin{subequations} \label{eq:online-optimization}
\begin{align} 
u^*_t \in \arg\min_{u \in \mathbb{U}} ~&~ \lVert u - \pi_{\text{nom}}(x_t) \rVert^2\\
\mathst ~&~ \mathcal{H} (x_t, u) \geqslant L_hw_{\text{max}} \label{eq:online-optimization-const}. 
\end{align}
\end{subequations} \\ [-0.95cm]

An important challenge is synthesizing an R-DTCBF. In the following, we explain a novel synthesis problem. 
\begin{problem}[Synthesis] \label{problem:synthesis}
    Consider the system \eqref{eq:discrete-time-dynamical-system} with $\mathbb{U}$, $\mathbb{W}$, and the safe set $\mathcal{S}$ in \eqref{eq:safe-set}. The objective is to synthesize a function $h:\mathbb{R}^n \rightarrow \mathbb{R}$ that is Lipschitz continuous on its zero-superlevel set $\mathcal{C}$ and a $\gamma \in \mathcal{K}_{\infty}^{\leqslant \text{id}}$ such that
    \begin{enumerate}
        \item $(h, \gamma)$ satisfies Lemma \ref{lemma:sufficient-RDTCBF}, thus it is an R-DTCBF pair. \label{condition:DCBF-constraint}
        \item The zero-superlevel set $\mathcal{C}$ of $h$ is a subset of $\mathcal{S}$. \label{condition:safety}
    \end{enumerate}
\end{problem}
To address Problem \ref{problem:synthesis}, we introduce a counterexample-guided approach using the following assumptions:
\begin{assumption} \label{assumption:safe-set}
    The safe set $\mathcal{S}$ is bounded. 
\end{assumption}
\begin{assumption} \label{assumption:diff}
    The mappings $f$ and $s_i$, $i \in \{1,2, \hdots, n_s\}$, associated with the system \eqref{eq:discrete-time-dynamical-system} and $\mathcal{S}$ in \eqref{eq:safe-set} are continuous.
\end{assumption}

\begin{assumption} \label{assumption:control-admissible}
    The control admissible set $\mathbb{U}$ is defined as \vspace{-0.2cm}
\setlength{\belowdisplayskip}{1pt}
    \begin{align} \label{eq:control-admissible-set}
        \mathbb{U} \hspace{-0.06cm} \defeq \hspace{-0.06cm} \{ u \in \mathbb{R}^m \hspace{-0.04cm} \mid \hspace{-0.04cm} u_{i, \text{min}} \hspace{-0.05cm} \leqslant \hspace{-0.05cm}  u_i \hspace{-0.05cm}\leqslant \hspace{-0.05cm}u_{i, \text{max}},\hspace{0.04cm}i \in \{1,2, \hdots, m\} \},
    \end{align}
    where $u_{i, \text{min}}, u_{i, \text{max}} \in \mathbb{R}$ such that $u_{i, \text{min}} \leqslant u_{i, \text{max}}$.
\end{assumption}
Assumptions \ref{assumption:safe-set} and \ref{assumption:diff} are required for the verification algorithm, and Assumption \ref{assumption:control-admissible} is required for both the verification and training algorithms. In Section \ref{sec:RDCBF-Verification}, we present our verification algorithm, and in \mbox{Section \ref{sec:RDCBF-Train}}, we propose our counterexample-guided approach that uses the verification algorithm as a subroutine.
\vspace{-0.12cm}
\section{Verification Algorithm} \label{sec:RDCBF-Verification}
The overall verification algorithm for a given pair $(h,\gamma)$, where $h$ is continuously differentiable, consists of two sub-algorithms: one for verifying the conditions of Lemma \ref{lemma:sufficient-RDTCBF} (\textit{verification of a candidate R-DTCBF}), and one for verifying $\mathcal{C} \subseteq \mathcal{S}$ (\textit{verification of safety}). 
For the verification of a candidate R-DTCBF, we extend the algorithm proposed in our previous work \cite{Shakhesi2024} to account for disturbances in the system~\eqref{eq:discrete-time-dynamical-system}. 
This sub-algorithm is introduced in \mbox{Section \ref{sec:verification-R-DTCBF}}, and the overall algorithm is proposed in Section \ref{sec:verification-overall}.
For the verification of safety, that is, verifying whether $\mathcal{C} \subseteq \mathcal{S}$, this is equivalent to checking that \mbox{$s_i^* \geqslant 0$} for all $i \in \{1,2,\hdots, n_s\}$, where \mbox{$s_i^* \defeq \min_{x \in \mathcal{C}} s_i(x)$} and the mappings $s_i$ define the safe set $\mathcal{S}$ in \eqref{eq:safe-set}. Since this verification can be performed using either the $\alpha$BB algorithm \cite{Floudas2010a} or the method in Section \ref{sec:verification-R-DTCBF}, we omit further details on the \texttt{Verify\_Safety} function in the overall algorithm (Algorithm \ref{algorithm:overall-verification}) due to space constraints.
\vspace{-0.1cm}

\vspace{-0.025cm}
\subsection{Verification of a Candidate R-DTCBF} \label{sec:verification-R-DTCBF}
In this approach, we first consider a hyperrectangle \vspace{-0.1cm}
\setlength{\belowdisplayskip}{3.5pt}
\begin{align}
    \mathbb{X} \defeq \{x \in \mathbb{R}^n \mid x_i^{lb} \leqslant x_i \leqslant x_i^{ub}, ~ i\in \{1,2,\hdots, n\}\}, \nonumber
\end{align} 
denoted by $\mathbb{X} \hspace{-0.08cm}\defeq \hspace{-0.08cm} [x^{lb},\hspace{0.05cm}x^{ub}]$, where $x^{lb} \hspace{-0.05cm},x^{ub} \hspace{-0.08cm}\in \hspace{-0.04cm} \mathbb{R}^n$ with $x_i^{lb} \hspace{-0.06cm} \leqslant \hspace{-0.06cm} x_i^{ub}\hspace{-0.08cm}$, $i \in \{1,2, \hdots, n\}$, such that $\mathcal{S} \subseteq \mathbb{X}$. 
We also assume that $\widetilde{L}_h \in \mathbb{R}_{\geqslant 0}$ is given as a \vspace{0.025cm}candidate Lipschitz constant for $h$ on the set $\mathcal{C}$. If $\widetilde{L}_h$ is not provided, refer to \mbox{Remark \ref{remark:not-given-Lip}}. 

We determine a control input that satisfies the \mbox{R-DTCBF} constraint \eqref{eq:DTCBF-constraint} at the middle of $\mathbb{X}$ by solving an optimization problem. If this control input satisfies \eqref{eq:DTCBF-constraint} for all $x \in \mathbb{X}$, and $\widetilde{L}_h$ is a Lipschitz constant for $h$ on $\mathbb{X}$, then $(h, \gamma)$ is verified as an R-DTCBF pair that satisfies Lemma \ref{lemma:sufficient-RDTCBF}. Otherwise, we divide $\mathbb{X}$ into smaller subdomains and repeat the same procedure within each one. 
To simplify the introduction of the algorithm (detailed in Algorithm \ref{algorithm:verification-RDCBF}), we initially focus on a three-step approach within \mbox{\textit{a subdomain}} $\mathbb{X}^{(k)} \defeq [x^{lb, (k)}, ~ x^{ub, (k)}] \subseteq \mathbb{X}$, where $k\in \mathbb{N}_0$ is the subdomain number.
\begin{enumerate}[label=\underline{Step \Roman{enumi}}:,ref=\Roman{enumi}, wide=0pt, listparindent=\parindent] 
\item \label{algorithm:verification:unknown:step1} We select the state $\bar{x}^{(k)} \in \mathbb{R}^n$ as the middle of $\mathbb{X}^{(k)}$, \vspace{-0.15cm} 
\begin{align} 
    \bar{x}^{(k)}_i \defeq \frac{x^{ub, (k)}_i + x^{lb, (k)}_i}{2}, 
\end{align} \\ [-0.2cm]
where $i \in \{1,2, \hdots, n\}$ and $\bar{x}^{(k)} \defeq [\bar{x}^{(k)}_1~ \hdots~ \bar{x}^{(k)}_n ]^{\T}$.
\item \label{algorithm:verification:unknown:step2} For $\bar{x}^{(k)}$, we aim to find an admissible control input $u^{*,(k)} \in \mathbb{U}$ that satisfies the \mbox{R-DTCBF constraint \eqref{eq:DTCBF-constraint}}. Therefore, we solve the following optimization problem to global optimality using the $\alpha$BB algorithm \cite{Floudas2010a}: \vspace{-0.1cm}
\begin{align} 
    u^{*, (k)} \in \underset{u \in \mathbb{U}}{\argmax} ~&~ \underbrace{\mathcal{H} (\bar{x}^{(k)},u) - \widetilde{L}_hw_{\text{max}} }_{\defeq \mathcal{F}_u(\bar{x}^{(k)},u)}. \label{eq:verification:unknown:step2-objective}
\end{align}
If $\mathcal{F}_u(\bar{x}^{(k)},u^{*, (k)}) < 0$ and $\bar{x}^{(k)} \in \mathcal{C}$, then $\bar{x}^{(k)}$ is a counterexample, and the algorithm terminates. Otherwise, we proceed to Step \ref{algorithm:verification:unknown:step3}.
\item \label{algorithm:verification:unknown:step3} 
A Lipschitz constant for $h$ on $\mathbb{X}^{(k)}$ can be obtained by solving $\max_{x \in \mathbb{X}^{(k)}}\lVert \nabla_x h(x) \rVert$. Therefore, to verify whether $\widetilde{L}_h$ is a Lipschitz constant for $h$ on $\mathbb{X}^{(k)}$, we check \mbox{$-\widebreve{\mathcal{F}}_1^{*, (k)} \leqslant \widetilde{L}_h$} by solving the convex optimization problem
\vspace{-0.05cm}
\begin{subequations} \label{eq:sec3:verification:unknown:step3-disturbance}
    \begin{align} 
        \widebreve{\mathcal{F}}_1^{*, (k)} \defeq \min_{x \in \mathbb{X}^{(k)}} & 
        \widebreve{\mathcal{F}}_1^{(k)}(x)\\
        \mathst ~&~ \widebreve{h}_{-}^{(k)}(x) \leqslant 0. 
    \end{align}
\end{subequations} \\ [-0.4cm]
Here, $\widebreve{\mathcal{F}}_1^{(k)}$ is a convex underestimator of $-\lVert \nabla_x h(x) \rVert$ within $\mathbb{X}^{(k)}$, i.e., $\widebreve{\mathcal{F}}_1^{(k)}(x) \leqslant -\lVert \nabla_x h(x) \rVert$ holds for all $x \in \mathbb{X}^{(k)}$ and $\widebreve{\mathcal{F}}_1^{(k)}$ is convex within $\mathbb{X}^{(k)}$, and is constructed as \vspace{-0.1cm}
\begin{align}
\widebreve{\mathcal{F}}_1^{(k)} &\hspace{-0.05cm} \defeq \hspace{-0.05cm} -\lVert \nabla_x h(x) \rVert \hspace{-0.05cm} + \hspace{-0.05cm} \sum_{i=1}^{n} \hspace{-0.05cm} \alpha_{\mathcal{F}_1,i}^{(k)}\bigl(x_i^{lb,(k)} \hspace{-0.1cm} - x_i\bigr)\bigl(x_i^{ub,(k)} \hspace{-0.1cm} - x_i\bigr), \nonumber
\end{align} \\ [-0.1cm]
where $\alpha_{\mathcal{F}_1,i}^{(k)} \in \mathbb{R}_{\geqslant 0}$, $i \hspace{-0.05cm} \in \hspace{-0.05cm} \{1,2, \hdots, n\}$, are sufficiently large and are obtained using the methods discussed in \cite{Floudas2010a}. Similarly, $\widebreve{h}_{-}^{(k)}$ is a convex underestimator of $-h$ within $\mathbb{X}^{(k)}$.

\begin{remark} \label{remark:not-given-Lip}
    If $\widetilde{L}_h$ is not provided, a Lipschitz constant for $h$ on $\mathcal{C}$ can be obtained and updated in each iteration by taking the maximum value of $-\widebreve{\mathcal{F}}^{*,(k)}_1$ over all subdomains. 
\end{remark}

We then check whether $u^{*,(k)}$ can satisfy the R-DTCBF constraint \eqref{eq:DTCBF-constraint} for all $x \in \mathbb{X}^{(k)}$. This involves verifying whether $\widebreve{\mathcal{F}}_2^{*,(k)} \geqslant 0$ by solving the convex optimization problem \vspace{-0.05cm}
\begin{subequations}\label{eq:verification:unknown:step3-convex}
    \begin{align} 
        \widebreve{\mathcal{F}}_2^{*, (k)}  \defeq \min_{x \in \mathbb{X}^{(k)}}  & \widebreve{\mathcal{F}}_2^{(k)}(x) \label{eq:verification:unknown:step3-convex-objective} \\
        \mathst ~& \widebreve{h}_{-}^{(k)}(x) \leqslant 0. \label{eq:verification:unknown:step3-convex-constraint}
    \end{align}
\end{subequations} \\ [-0.4cm]
Here, $\widebreve{\mathcal{F}}_2^{(k)}$ is a convex underestimator of 
    $\mathcal{H}(x,u^{*,(k)}) - \widetilde{L}_hw_{\text{max}}$ 
within $\mathbb{X}^{(k)}$. 
\end{enumerate}

By solving \eqref{eq:sec3:verification:unknown:step3-disturbance} and \eqref{eq:verification:unknown:step3-convex} within a subdomain $\mathbb{X}^{(k)}$, we encounter three cases:
\begin{enumerate}[label=(\Alph*), ref=(\Alph*), leftmargin=18pt]
    \item  $\widetilde{L}_h \geqslant -\widebreve{\mathcal{F}}_1^{*, (k)}$ and $\widebreve{\mathcal{F}}_2^{*, (k)} \geqslant 0$: $(h, \gamma)$ is verified as an R-DTCBF pair within the subdomain $\mathbb{X}^{(k)}$. Thus, we approve $\mathbb{X}^{(k)}$ (Line \ref{alg:line-17} of Algorithm \ref{algorithm:verification-RDCBF} below). \label{case:convex-verification-known-A} 
    \item Problem \eqref{eq:sec3:verification:unknown:step3-disturbance} is infeasible: $\mathbb{X}^{(k)}$ is entirely outside of the zero-superlevel set $\mathcal{C}$ of $h$, as $\widebreve{h}^{(k)}_{-}(x) \leqslant -h(x)$ for all $x \in \mathbb{X}^{(k)}$. Thus, we approve $\mathbb{X}^{(k)}$  (Line \ref{alg:line-17}). \label{case:convex-verification-known-B}
    \item \label{case:convex-verification-known-C} 
    $ \widetilde{L}_h < -\widebreve{\mathcal{F}}_1^{*, (k)}$ or $\widebreve{\mathcal{F}}_2^{*, (k)} < 0$: We divide $\mathbb{X}^{(k)}$ into smaller subdomains (Lines \ref{alg-line-21} and \ref{alg-line-22}). 
\end{enumerate}
We continue this process until either all subdomains of $\mathbb{X}$ are approved, showing that $(h, \gamma)$ is an R-DTCBF pair that satisfies Lemma \ref{lemma:sufficient-RDTCBF}, or a counterexample is found showing that $(h, \gamma)$ does not satisfy the conditions of Lemma \ref{lemma:sufficient-RDTCBF}. However, as discussed in \cite{Shakhesi2024}, to guarantee that the algorithm terminates in finite time, we must impose stopping criteria. Here, for simplicity, we only consider the size of the subdomains; for alternative stopping criteria, see \cite{Shakhesi2024}. Specifically, the algorithm terminates on $\mathbb{X}^{(k)}$ and provides $\bar{x}^{(k)}$ as a \textit{potential} counterexample for the R-DTCBF \mbox{constraint \eqref{eq:DTCBF-constraint}}, if \vspace{-0.1cm}
\begin{align}
    \sum_{i=1}^{n}\bigl(x_i^{ub,(k)} - x_i^{lb,(k)}\bigr)^2 \leqslant \epsilon, \label{eq:verification:stopping-criterion}
\end{align}
where $\epsilon \in \mathbb{R}_{> 0}$ is a predefined value.
\begin{remark}
   If \eqref{eq:verification:stopping-criterion} is satisfied, the algorithm may incorrectly falsify a valid R-DTCBF (though it never verifies a non-valid one). This does not affect the overall counterexample-guided approach in Section \ref{sec:counter-approach}, and $\epsilon$ can be reduced in the next iteration of the algorithm to achieve more accurate verification.
\end{remark}
\begin{algorithm}[H]
{\small
\caption{Algorithm for the Verification of a Candidate R-DTCBF (\texttt{Verify\_R-DTCBF})} \label{algorithm:verification-RDCBF}
\begin{algorithmic}[1]
\State \textbf{Input:} $f, \mathbb{U}, \mathbb{W}, h, \gamma, \epsilon, \widetilde{L}_h$. $\mathbb{X}$ such that \mbox{$\mathcal{S}\subseteq \mathbb{X}$}.
\State \textbf{Output:} $x_{\text{ce}, \text{R-DTCBF}}$ \Comment{$x_{\text{ce}, \text{R-DTCBF}}$ is a counterexample}
\State $\mathbb{X}^{(0)} \gets \mathbb{X}$
\State $\mathcal{D} \gets \{0\}$ \Comment{$\mathcal{D}$ is the list of subdomain numbers that remain to be handled}
\State $n_{\text{dom}} \gets 0$ \Comment{Number of subdomains}
\While{$\mathcal{D} \neq \emptyset $}
\State $k \gets \text{any subdomain number from the list } \mathcal{D}$ 
\State $\bar{x}^{(k)} \gets$ middle of $\mathbb{X}^{(k)}$
\State $u^{*,(k)} \gets$ solve \eqref{eq:verification:unknown:step2-objective} using $\bar{x}^{(k)}$
\If{$\mathcal{F}_u(\bar{x}^{(k)},u^{*,(k)}) < 0$ and $\bar{x}^{(k)} \in \mathcal{C}$}
    \State \textbf{return} $\bar{x}^{(k)}$
\EndIf
\State $\widebreve{\mathcal{F}}_1^{*,(k)}  \gets$ solve \eqref{eq:sec3:verification:unknown:step3-disturbance} within $\mathbb{X}^{(k)}$
\State $\widebreve{\mathcal{F}}_2^{*,(k)} \gets $ solve \eqref{eq:verification:unknown:step3-convex} within $\mathbb{X}^{(k)}$ using $u^{*,(k)}$
\If{Case \ref{case:convex-verification-known-A} or Case \ref{case:convex-verification-known-B}} 
    \State  $\mathcal{D} \gets \mathcal{D} \setminus \{k\}$ \label{alg:line-17}
\Else \Comment{Case \ref{case:convex-verification-known-C}}
\If{stopping criterion \eqref{eq:verification:stopping-criterion} is met for $\mathbb{X}^{(k)}$}
    \State \textbf{return} $\bar{x}^{(k)}$
\EndIf
\State divide $\mathbb{X}^{(k)}$ into $\mathbb{X}^{(n_{\text{dom}}+1)}$ and $\mathbb{X}^{(n_{\text{dom}}+2)}$ \label{alg-line-21}
\State $\mathcal{D} \gets (\mathcal{D} \setminus \{k\}) \cup \{n_{\text{dom}}+1,~ n_{\text{dom}}+2\}$ \label{alg-line-22}
\State $n_{\text{dom}} \gets n_{\text{dom}}+2$
\EndIf
\EndWhile
\State \textbf{return} $[~]$
\end{algorithmic}
}
\end{algorithm}

\subsection{Overall Verification Algorithm} \label{sec:verification-overall}
The overall algorithm first verifies safety, $\mathcal{C} \subseteq \mathcal{S}$ (\mbox{Line \ref{alg:safety-Line3}} in Algorithm \ref{algorithm:overall-verification}). If a counterexample $x_{\text{ce,safety}}$ is found, it is returned, and the algorithm terminates. Otherwise, it proceeds to verify the conditions of Lemma \ref{lemma:sufficient-RDTCBF} using \mbox{Algorithm \ref{algorithm:verification-RDCBF}} (Line~\ref{alg:safety-Line7} in \mbox{Algorithm \ref{algorithm:overall-verification}}). If a counterexample $x_{\text{ce,R-DTCBF}}$ is found, it is returned, and the algorithm terminates. Otherwise, $(h, \gamma)$ is an R-DTCBF pair that satisfies Lemma \ref{lemma:sufficient-RDTCBF} and $\mathcal{C} \subseteq \mathcal{S}$.
\begin{algorithm}
{\small
\caption{Overall Verification Algorithm (\texttt{Verify})} \label{algorithm:overall-verification}
\begin{algorithmic}[1] 
\State \textbf{Input:} $f,\mathbb{U}, \mathbb{W}, h, \gamma,\epsilon,\widetilde{L}_h,\mathcal{S}$. $\mathbb{X}$ such that $\mathcal{S}\subseteq \mathbb{X}$.
\State \textbf{Output:} $\texttt{verif}, x_{\text{ce}, \text{safety}}, x_{\text{ce}, \text{R-DTCBF}}$ \Comment{If $(h,\gamma)$ is verified, $\texttt{verif}= \text{True}$. Otherwise, \texttt{verif} = False.}
\State $x_{\text{ce}, \text{safety}} \gets \texttt{Verify\_Safety}(h, \mathcal{S},\mathbb{X})$ \label{alg:safety-Line3}
\If{$x_{\text{ce}, \text{safety}}$ is non-empty}
    \State \textbf{return} $\text{False}, x_{\text{ce}, \text{safety}}, [~]$
\Else
    \State $x_{\text{ce}, \text{R-DTCBF}} \hspace{-0.05cm}\gets \hspace{-0.05cm}\texttt{Verify\_R-DTCBF}(f, \mathbb{U}, \mathbb{W}, h, \gamma, \epsilon, \widetilde{L}_h, \mathbb{X})$ \label{alg:safety-Line7}
    \If{$x_{\text{ce}, \text{R-DTCBF}}$ is non-empty}
        \State \textbf{return} $\text{False}, [~], x_{\text{ce}, \text{R-DTCBF}}$
    \Else 
        \State \textbf{return} $\text{True}, [~], [~]$
    \EndIf
\EndIf
\end{algorithmic}
}
\end{algorithm}

\section{Training R-DTCBFs} \label{sec:RDCBF-Train}
First, in Section \ref{sec:training-candidate}, we propose a method to train a candidate R-DTCBF pair $(h,\gamma)$ with a control policy $\pi$ as a friend such that both the conditions of Lemma \ref{lemma:sufficient-RDTCBF} and $\mathcal{C} \subseteq \mathcal{S}$ hold on a given set of sampled states. Then, in Section \ref{sec:counter-approach}, we propose our counterexample-guided approach to learn an R-DTCBF that is verified using the verification algorithm.
\subsection{Training Candidate R-DTCBFs} \label{sec:training-candidate}
\subsubsection{Parameterization}
We parameterize $h(x;\vartheta)$ using a set of continuously differentiable basis functions and the vector of unknown coefficients \mbox{$\vartheta \in \mathbb{R}^{n_\vartheta}$}. 
Note that keeping the parameterization of $h$ simple makes the online optimization problem \eqref{eq:online-optimization} for synthesizing a safe controller more practical by ensuring reasonable computational time. Therefore, we opt not to parameterize $h$ as a neural network (NN).

In contrast, the control policy $\pi$ is only required to train an R-DTCBF and is not needed afterward. As a result, similar to \cite{Formally-Verified-NN} and \cite{Htanh}, we parameterize the control \mbox{policy $\pi_i$}, \mbox{$i \in \{1,2, \hdots, m\}$}, with $\pi \defeq [\pi_1 ~ \hdots ~ \pi_m]^{\T}$ as the NN \vspace{-0.05cm} \setlength{\belowdisplayskip}{6pt}
\[
\left\{
\begin{aligned}
    &y^{(0)} \defeq x, \\
    &y^{(j)} \defeq \mathrm{ReLU}(W^{(j)}y^{(j-1)} + b^{(j)}), ~ j \in \{1,2, \dots, l-1\}, \\
    &\pi_i(x; \kappa) \defeq \mathrm{Hardtanh} (W^{(l)} y^{(l-1)} + b^{(l)}, u_{i, \text{min}}, u_{i, \text{max}}),
\end{aligned}
\right.
\]
where $\kappa \defeq \{W^{(j)}, b^{(j)}\}_{j=1}^{l}$ represents the set of weight matrices and biases, and $l \in \mathbb{N}$ is the number of layers. Here, \mbox{$W^{(j)} \in \mathbb{R}^{d_{j} \times d_{j-1}}$} is a weight matrix, $b^{(j)} \in \mathbb{R}^{d_{j}}$ is a bias vector, and $d_j \in \mathbb{N}$ is the number of neurons in the $j$-th layer, $j \in \{0,1, \hdots, l\}$. The number of neurons in the input and output layers are $d_0 \defeq n$ and $d_l \defeq 1$, respectively, while the number of neurons in the hidden layers (i.e., for \mbox{$j \in \{1,2, \hdots, l-1\}$}) are user-specified. Moreover, the activation function $\mathrm{ReLU}$ (rectified linear unit), defined as $\mathrm{ReLU}(y) \defeq \max(0,y)$, is applied element-wise in the hidden layers. For the output layer, we use $\mathrm{Hardtanh}$ as the activation function, defined as \vspace{-0.15cm}
\begin{align}
\mathrm{Hardtanh}(u_i, u_{i,\text{min}}, u_{i,\text{max}}) \hspace{-0.07cm} \defeq \hspace{-0.07cm}\left\{
\begin{aligned}
    &u_{i,\text{max}}, \hspace{0.75cm}u_{i, \text{max}} \leqslant u_i, \\
    &u_i, \hspace{0.51cm}u_{i, \text{min}} \leqslant u_i \leqslant u_{i, \text{max}}, \\
    &u_{i, \text{min}}, \hspace{1.298cm}u_i \leqslant u_{i, \text{min}}, 
\end{aligned}
\right. \nonumber
\end{align}
to ensure that the control policy $\pi$ is admissible on the sampled states with respect to $\mathbb{U}$ in \eqref{eq:control-admissible-set}. 
Regarding \mbox{$\gamma \in \mathcal{K}^{\leqslant \text{id}}_{\infty}$}, we parameterize $\gamma$ as a linear function for simplicity, i.e., \vspace{-0.1cm}
\begin{align} \label{eq:parameterized-gamma}
    \gamma(r) \defeq \gamma_0r, \quad r \in \mathbb{R}_{\geqslant 0},
\end{align}
where $\gamma_0 \in \mathbb{R}_{> 0}$ is an unknown constant that satisfies \mbox{$\gamma_{0,\text{min}} \leqslant \gamma_0 \leqslant \gamma_{0,\text{max}} \leqslant 1$}, with $\gamma_{0,\text{min}}, \gamma_{0,\text{max}} \in \mathbb{R}_{>0}$ being user-specified constants.

\subsubsection{Sampling and Designing the Loss Function}
We require a finite set of sampled states \mbox{$\mathcal{X}_s \subset \mathcal{S}$} for the \mbox{conditions} of Lemma \ref{lemma:sufficient-RDTCBF} and also a finite set of sampled states \mbox{$\mathcal{X}_u \subset \mathbb{R}^n \setminus \mathcal{S}$} for $\mathcal{C} \subseteq \mathcal{S}$. These sampled states are \mbox{initially} selected randomly within their respective domains, after which counterexamples are iteratively added. In the \mbox{following}, we introduce loss terms that, when minimized to global optimality, promote satisfaction of the conditions in Lemma \ref{lemma:sufficient-RDTCBF} and \mbox{$\mathcal{C}\subseteq \mathcal{S}$} on the sampled states. Ideally, the global minimum of each loss term should be zero to prevent trade-offs between their satisfaction. Therefore, we use the $\mathrm{ReLU}$ function in each loss term.

Firstly, to promote the Lipschitz continuity of the trained function $h$ on the set $\mathcal{C}$, we introduce the loss term \vspace{-0.1cm}
\begin{align}
    &\mathcal{L}_1(\mathcal{X}_s; \vartheta) \defeq \nonumber \\
    &\sum_{x \in \mathcal{X}_s} \mathrm{ReLU}\left(\mathrm{ReLU}\bigl(h(x;\vartheta) + \eta\bigr) \bigl(\lVert \nabla_x h(x;\vartheta)\rVert - \widetilde{L}_h\bigr) \right) \nonumber
\end{align}
to enforce $\lVert \nabla_x h(x;\vartheta) \rVert \leqslant \widetilde{L}_h$ for all $x \in \mathcal{X}_s$ satisfying $h(x;\vartheta) > -\eta$, where $\eta \in \mathbb{R}_{>0}$ is sufficiently small. Here, \mbox{$\widetilde{L}_h \in \mathbb{R}_{\geqslant 0}$} is a predefined value that serves as a Lipschitz constant for the resulting $h$ on its zero-superlevel \mbox{set $\mathcal{C}$} and is estimated based on the basis functions of $h$. Regarding Lemma~\ref{lemma:sufficient-RDTCBF}, since the R-DTCBF \mbox{constraint \eqref{eq:DTCBF-constraint}} only needs to hold within $\mathcal{C}$, which is unknown, we design the loss term \vspace{-0.1cm}
\begin{align}
    &\mathcal{L}_2(\mathcal{X}_s; \vartheta, \gamma_0, \kappa) \defeq  \nonumber \\
    &\hspace{-0.1cm} \sum_{x \in \mathcal{X}_s} \hspace{-0.1cm} \mathrm{ReLU}\left(\Gamma \bigl(h(x; \vartheta) \bigr)\left(\mathcal{H}\bigl(x,\pi(x;\kappa); \vartheta, \gamma_0\bigr) - \widetilde{L}_hw_{\text{max}} \right)\right), \nonumber
\end{align}
where \vspace{-0.17cm}
\begin{align}
   \Gamma(z) \defeq \frac{-1}{1+ e^{-c_1z}} + c_2. \nonumber
\end{align}
Here, $c_1, c_2 \in \mathbb{R}_{>0}$ are predetermined constants such that $h(x;\vartheta) > -\zeta$ implies $\Gamma(h(x;\vartheta)) < 0$, and \mbox{$h(x;\vartheta) \leqslant -\zeta$} implies $\Gamma(h(x;\vartheta)) \geqslant 0$, where $\zeta \in \mathbb{R}_{>0}$ is adjusted by the values of $c_1$ and $c_2$ and chosen to be sufficiently small. Therefore, if a sampled state is within the set \mbox{$\mathcal{C}_\zeta \defeq\{x \in \mathbb{R}^n \mid h(x; \vartheta) > -\zeta\}$}, then the R-DTCBF \mbox{constraint \eqref{eq:DTCBF-constraint}} must be satisfied to prevent an increase in the loss term. Conversely, if a sampled state lies outside $\mathcal{C}_\zeta$, \eqref{eq:DTCBF-constraint} must not be satisfied. In this case, the state does not enter the set $\mathcal{C}$ at the next time instant. This steers the training of a candidate R-DTCBF to a \mbox{larger set $\mathcal{C}$}. A larger $\mathcal{C}$ reduces the conservatism of the safe controller synthesis (i.e., optimization problem \eqref{eq:online-optimization}) and improves overall performance, for instance reducing tracking errors. Moreover, to ensure that $\mathcal{C}$ is sufficiently large, we incorporate an additional loss term. If computing the area $\phi(\vartheta)$ of the zero-superlevel set of the parameterized $h$ is straightforward, we define \vspace{-0.09cm} \setlength{\belowdisplayskip}{4.8pt}
\begin{align}
    \mathcal{L}_3(\vartheta) \defeq \mathrm{ReLU}\bigl(-\phi(\vartheta) + \phi_0\bigr), \nonumber
\end{align}
to enforce $\phi(\vartheta) \geqslant \phi_0$, where $\phi_0 \in \mathbb{R}_{>0}$ is a predefined value. If computing $\phi(\vartheta)$ is challenging, an alternative approach is to consider a finite set $\mathcal{X}_{I}$ of initial states such that the resulting R-DTCBF $h$ remains non-negative for all states in this set. In this case, we consider the loss term \vspace{-0.08cm}
\begin{align}
    \mathcal{L}_3(\mathcal{X}_I; \vartheta) \defeq \sum_{x \in \mathcal{X}_I} \mathrm{ReLU}\bigl(-h(x; \vartheta)\bigr). \nonumber
\end{align}

Moreover, to ensure that $\mathcal{C} \subseteq \mathcal{S}$ holds on the sampled states, we consider the loss term \vspace{-0.08cm}
\begin{align}
    \mathcal{L}_4(\mathcal{X}_u; \vartheta) &\defeq \sum_{x \in \mathcal{X}_u} \mathrm{ReLU}\bigl(h(x; \vartheta) + \delta \bigr), \nonumber
\end{align}
where $\delta \in \mathbb{R}_{>0}$ is chosen sufficiently small.
Lastly, regarding $\gamma$ parameterized as in \eqref{eq:parameterized-gamma}, we need to enforce that the unknown coefficient $\gamma_0$ satisfies \mbox{$\gamma_{0,\text{min}} \leqslant \gamma_0 \leqslant \gamma_{0,\text{max}} \leqslant 1$}. Thus, we consider the loss term
\begin{align}
    \mathcal{L}_5(\gamma_0) \defeq \mathrm{ReLU}\bigl(-\gamma_0 + \gamma_{0,\text{min}}\bigr) + \mathrm{ReLU}\bigl(\gamma_0 - \gamma_{0,\text{max}}\bigr). \nonumber
\end{align}
The loss function is then a weighted sum of the above-mentioned loss terms, \vspace{-0.15cm}
\begin{align}
    &\mathcal{L}(\mathcal{X}_s, \mathcal{X}_u; \vartheta, \gamma_0,\kappa) \defeq \alpha_1 \mathcal{L}_1(\mathcal{X}_s; \vartheta) + \alpha_2 \mathcal{L}_2(\mathcal{X}_s; \vartheta, \gamma_0, \kappa) + \nonumber\\
    &\qquad \qquad \qquad \qquad \alpha_3 \mathcal{L}_3(\vartheta) + \alpha_4 \mathcal{L}_4(\mathcal{X}_u; \vartheta) + \alpha_5 \mathcal{L}_5(\gamma_0), \nonumber
\end{align}
where $\alpha_1, \hdots, \alpha_5 \in \mathbb{R}_{> 0}$ are predefined.
\begin{theorem}
    If $\mathcal{L}(\mathcal{X}_s, \mathcal{X}_u;\vartheta, \gamma_0, \kappa) = 0$, the trained pair $(h(x;\vartheta), \gamma)$ satisfies the conditions of Lemma \ref{lemma:sufficient-RDTCBF} and $\lVert \nabla_x h(x;\vartheta) \rVert \leqslant \widetilde{L}_h$ for all \mbox{$x \in \mathcal{X}_s$} satisfying \mbox{$h(x;\vartheta) \geqslant 0$}, and $h(x; \vartheta) < 0$ holds for all $x \in \mathcal{X}_u$. \\
    \begin{proof}
        It follows directly from the construction of $\mathcal{L}$.
    \end{proof}
\end{theorem}

The training procedure for the candidate R-DTCBF pair $(h,\gamma)$ and its friend $\pi$ is outlined in Algorithm \ref{algorithm:learning-candidate}. Following existing works, the sampled states are shuffled (Line \ref{alg:learning-line-7}), split into batches (Line \ref{alg:learning-line-8}), and used to compute the gradient of the loss function and update the parameters (Line \ref{alg:learning-line-14}) until either the loss function converges to zero for all sampled states or \texttt{max\_epoch} is reached. If the latter occurs, the parameters are reinitialized and training is repeated.
\begin{algorithm}[H]
{\small
\caption{Training Candidate R-DTCBFs (\texttt{Train})} \label{algorithm:learning-candidate}
\begin{algorithmic}[1] 
\State \textbf{Input:} $f, \mathbb{U}, \mathbb{W}, \mathcal{X}_s,\mathcal{X}_u$, \texttt{batch\_size}, \texttt{max\_epoch}
\State \textbf{Output:} Candidate $h$, $\gamma$
\State {\footnotesize $\texttt{n\_batches} \hspace{-0.09cm}\gets \hspace{-0.09cm} \max(\text{len}(\mathcal{X}_s) / \texttt{batch\_size}, \text{len}(\mathcal{X}_u) / \texttt{batch\_size})$} \Comment{len represents the number of sampled states}
\While {true} \Comment{Infinite Loop}
    \State $(\vartheta^{(0)}, \gamma_0^{(0)}, \kappa^{(0)}) \gets \texttt{Initialize()}$ \label{alg:learning-line-5} 
    \State $k \gets 0$ 
    \For{$\texttt{epoch} \gets 1$ \textbf{ to } $\texttt{max\_epoch}$}
        \State Shuffle $\mathcal{X}_s$ and $\mathcal{X}_u$ \label{alg:learning-line-7}
        \State Split $\mathcal{X}_s$ and $\mathcal{X}_u$ into batches of size $\texttt{batch\_size}$ \label{alg:learning-line-8}
        \For{$i \gets 1$ \textbf{ to } $\texttt{n\_batches}$}
            \State Get batch $\mathcal{B}^{(i)}_s$ and $\mathcal{B}^{(i)}_u$ from $\mathcal{X}_s$ and $\mathcal{X}_u$, resp.
            \State $(\vartheta^{(k+1)}, \gamma_0^{(k+1)}, \kappa^{(k+1)}) \gets (\vartheta^{(k)}, \gamma_0^{(k)}, \kappa^{(k)}) - \mu \nabla_{\vartheta,\gamma_0,\kappa} \mathcal{L}\bigl(\mathcal{B}^{(i)}_s, \mathcal{B}^{(i)}_u; \vartheta, \gamma_0, \kappa \bigr)\big|_{\vartheta=\vartheta^{(k)}, \gamma_0=\gamma_0^{(k)}, \kappa=\kappa^{(k)}}$ \label{alg:learning-line-14}  
            \State $k \gets k+1$
        \EndFor
        \If{$\mathcal{L}(\mathcal{X}_s, \mathcal{X}_u; \vartheta^{(k)}, \gamma_0^{(k)}, \kappa^{(k)}) = 0$}
            \State \textbf{return} $h(x; \vartheta^{(k)}), \gamma^{(k)}_0$
        \EndIf
    \EndFor
\EndWhile
\end{algorithmic}
}
\end{algorithm}

\subsection{Counterexample-Guided Approach} \label{sec:counter-approach}
The counterexample-guided approach, as detailed in Algorithm \ref{algorithm:counter-example}, employs the proposed verification algorithm iteratively. If a counterexample is generated that violates safety ($\mathcal{C} \subseteq \mathcal{S}$) or violates the conditions of Lemma \ref{lemma:sufficient-RDTCBF} (within a predefined tolerance), it is added to $\mathcal{X}_u$ and $\mathcal{X}_s$, respectively, for retraining. This process repeats until an R-DTCBF pair $(h,\gamma)$ is obtained that satisfies Lemma \ref{lemma:sufficient-RDTCBF} and $\mathcal{C} \subseteq \mathcal{S}$ for the system \eqref{eq:discrete-time-dynamical-system} with $\mathbb{U}$ and $\mathbb{W}$, or until a predefined maximum time is reached, indicating that the approach is not successful. 
\begin{remark}
Note that a trained control \mbox{policy $\pi$} is not used to verify $(h, \gamma)$ because training a control policy that is valid for all states would require a significantly larger number of sampled states, which we aim to avoid. Moreover, using $\pi$ for verification would be slower due to its potential complexity.
\end{remark}
\begin{remark}
    When retraining $h$, $\gamma$, and $\pi$ with Algorithm~\ref{algorithm:learning-candidate}, the parameters $\vartheta$, $\gamma_0$, and $\kappa$ are initialized with the values obtained from the previous counterexample-guided iteration.
\end{remark}
\setlength{\textfloatsep}{12pt} 
\begin{algorithm}
{\small
\caption{Counterexample-Guided Approach} \label{algorithm:counter-example}
\begin{algorithmic}[1] 
\State \textbf{Input:} $f, \mathbb{U}, \mathbb{W}, \mathcal{S}$, $\texttt{batch\_size}$, $\texttt{max\_epoch}$
\State \textbf{Output:} Verified $h, \gamma$
\State $\mathcal{X}_s, \mathcal{X}_u \gets \texttt{Initialize()}$
\State $\texttt{verif} \gets \text{False}$
\While {$\texttt{verif} = \text{False}$}
    \State $(h, \gamma)  \hspace{-0.07cm} \gets  \hspace{-0.12cm}$ \texttt{Train}{\footnotesize $(f,\mathbb{U},\mathbb{W},\mathcal{X}_s,\mathcal{X}_u,\texttt{batch\_size}, \texttt{max\_epoch})$}
    \State $\texttt{verif}, x_{\text{ce}, \text{safety}}, x_{\text{ce}, \text{R-DTCBF}} \hspace{-0.05cm}\gets \hspace{-0.05cm} \texttt{Verify}(f,\mathbb{U}, \mathbb{W},h, \gamma, \mathcal{S})$
    \State $\mathcal{X}_u \gets \mathcal{X}_u \cup \{x_{\text{ce}, \text{safety}}\}$
    \State $\mathcal{X}_s \gets \mathcal{X}_s \cup \{x_{\text{ce}, \text{R-DTCBF}}\}$
\EndWhile
\State \textbf{return} $h, \gamma$
\end{algorithmic}
}
\end{algorithm}

\section{Numerical Case Studies} \label{sec:numerical-case-studies}
The numerical case studies are solved using Python 3.12 and PyTorch 2.5 on a laptop with Linux Ubuntu 22.04, an Intel i7-12700H, and an Nvidia RTX A1000.
\subsection{Polynomial System Subject to Disturbances} 
\begin{figure}[ht]\centering 
\includegraphics[width=0.33\textwidth]{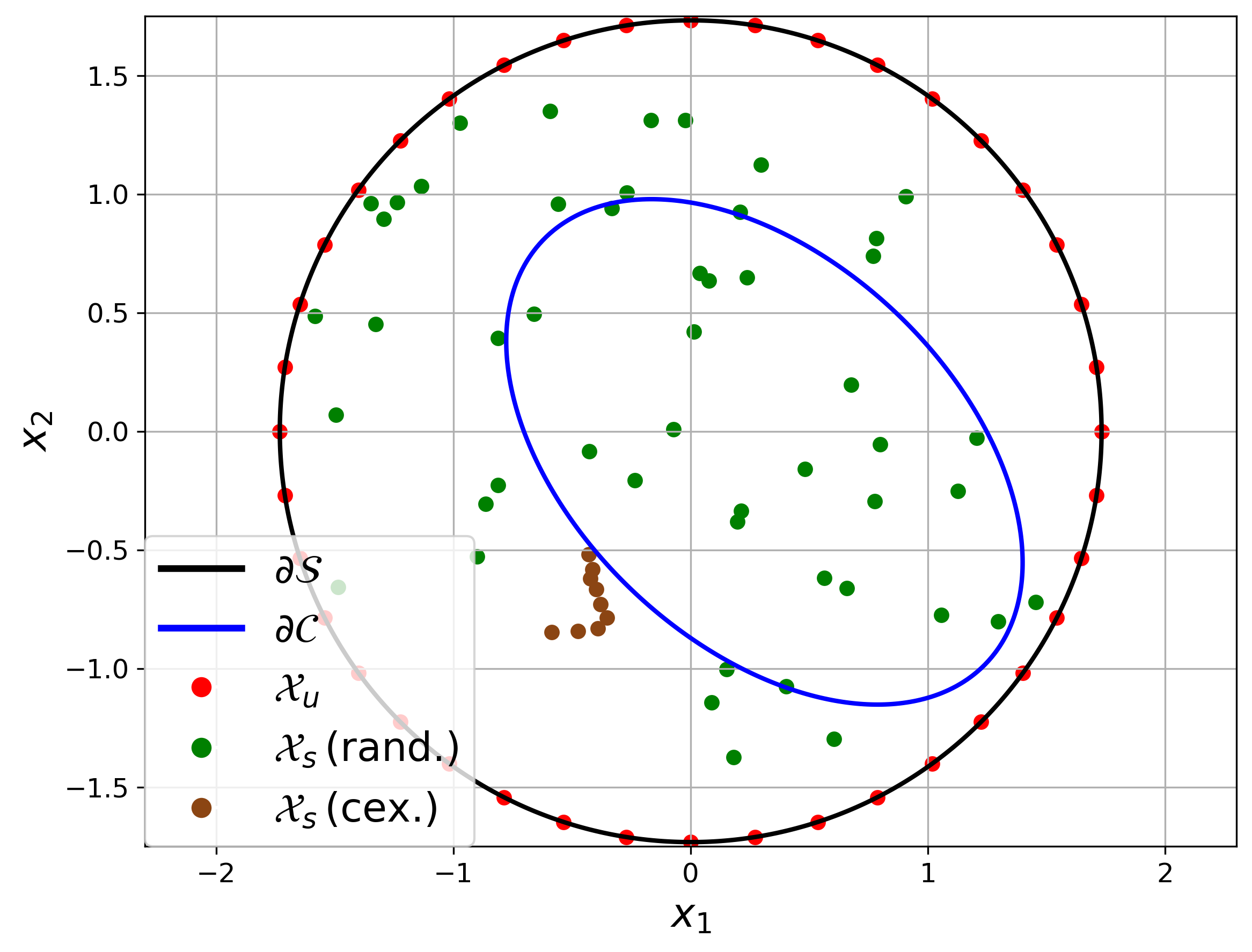}
    \vspace{-0.3cm}\caption{The black circle denotes the boundary of the safe set $\mathcal{S}$, and the blue ellipse represents the boundary of the zero-superlevel set of the R-DTCBF $h$. The red dots are sampled states in $\mathcal{X}_u$. The green and brown dots are sampled states that both belong to $\mathcal{X}_s$, representing randomly selected states and counterexamples, respectively.}
    \label{fig:polynomial-exm}
\end{figure}

We discretize the continuous-time polynomial system discussed in \cite{Wang2023a} using the forward Euler method with sampling time \mbox{$T_s \defeq \unit[0.1]{s}$} and consider additive bounded disturbances $w \defeq [w_1 ~ w_2]^{\T} \in \mathbb{R}^2$, \mbox{$\lVert w \rVert \leqslant 0.04$}, leading to 
\setlength{\belowdisplayskip}{4pt}
\begin{align}
    &\begin{bmatrix}
        x_1^+ \\
        x_2^+
    \end{bmatrix} = \begin{bmatrix}
        x_1 + x_2T_s \\
        x_2 + \bigl(x_1 + \frac{1}{3}x_1^3 + x_2\bigr)T_s
    \end{bmatrix} + \nonumber \\
    &\begin{bmatrix}
        \bigl(x_1^2 + x_2 + 1\bigr)T_s & 0 \\
        0 & \bigl(x_2^2 + x_1 + 1\bigr)T_s
    \end{bmatrix} \begin{bmatrix}
        u_1 \\
        u_2
    \end{bmatrix} + \begin{bmatrix}
        w_1 \\ 
        w_2
    \end{bmatrix}, \nonumber 
\end{align}
with the control input vector $u \defeq [u_1 ~u_2]^{\T} \in \mathbb{U}$ and 
\begin{align}
   \mathbb{U} \defeq \{ u \in \mathbb{R}^2 \mid -1.5 \leqslant u_1 \leqslant 1.5, ~ -1.5 \leqslant u_2 \leqslant 1.5 \}.
    \nonumber
\end{align}
The safe set $\mathcal{S}$ is defined as 
    $\mathcal{S} \defeq \left\{ x \in \mathbb{R}^2 \bigm| x_1^2 + x_2^2 \leqslant 3 \right \}$.
By parameterizing $h$ as a polynomial of degree two and $\gamma$ as in \eqref{eq:parameterized-gamma}, with $\gamma_{0,\text{min}} \defeq 0.7$ and $\gamma_{0,\text{max}} \defeq 0.9$, the \mbox{R-DTCBF} pair $(h, \gamma)$, depicted in Figure \ref{fig:polynomial-exm}, is learned using the proposed method after approximately $\unit[8]{min}$, as \vspace{-0.1cm}
\begin{align}
    h(x) \hspace{-0.05cm}&=\hspace{-0.05cm} -1.14x_1^2\hspace{-0.05cm} -\hspace{-0.05cm}1.02x_1x_2\hspace{-0.05cm} -\hspace{-0.05cm}1.19x_2^2\hspace{-0.05cm} +\hspace{-0.05cm} 0.62x_1\hspace{-0.05cm} +\hspace{-0.05cm} 0.11x_2\hspace{-0.05cm} +\hspace{-0.05cm} 1, \nonumber \\
    \gamma(r) &= 0.9r, ~~ r \in \mathbb{R}_{\geqslant 0}. \nonumber
\end{align}
\vspace{-0.7cm}
\vspace{-0.6cm}
\subsection{Cart-Pole}
\vspace{-0.05cm}
\hspace{-0.3cm}Consider the discretized cart-pole system discussed in \cite{Neural-Lyapunov} as 
\begin{align} 
    \begin{bmatrix}
        x_c^+ \\
        v_c^+ \\
        \theta^+ \\
        \omega^+
    \end{bmatrix} \hspace{-0.1cm}= \hspace{-0.1cm}\begin{bmatrix}
        x_c + v_cT_s \\
        v_c + \frac{T_s\left(u + m_p \sin\theta \hspace{0.02cm}(L\omega^2 - g\cos \theta)\right)}{m_c + m_p\sin^2 \theta} \\
        \theta + \omega T_s \\
        \omega + \frac{T_s\left(-u\cos \theta - m_pL\omega^2 \cos\theta \hspace{0.02cm} \sin \theta + (m_c + m_p)g \sin \theta\right)}{L(m_c + m_p\sin^2 \theta)}
    \end{bmatrix}, \nonumber
\end{align}
where $x_c, v_c \in \mathbb{R}$ are the horizontal position and velocity of the cart, respectively, $\theta \in \mathbb{R}$ is the angle of the pole from the upward vertical direction, and $\omega \in \mathbb{R}$ is its angular velocity. Moreover, $x \defeq [x_c ~ v_c ~ \theta ~ \omega]^{\T}$, and \mbox{$u \in \mathbb{U} \defeq [-30, ~ 30]$} is the horizontal force applied to the cart, which is considered as the control input to the system. We consider the following parameters: the cart mass $m_c \defeq \unit[2]{kg}$, the pole mass \mbox{$m_p \defeq \unit[0.1]{kg}$}, the pole length $L \defeq \unit[1]{m}$, and the sampling time $T_s \defeq \unit[0.01]{s}$.
The safe set $\mathcal{S}$ is defined as \vspace{-0.05cm}
\begin{align}
    \mathcal{S} \defeq \{ x \in \mathbb{R}^4 \bigm|\theta^2 + \omega^2 \leqslant (\pi/4)^2 \}. \nonumber
\end{align} \\ [-0.5cm]
By parameterizing $h$ as a polynomial of degree two and $\gamma$ as in \eqref{eq:parameterized-gamma},  the pair $(h,\gamma)$ is learned after about \unit[12]{min}, as \vspace{-0.05cm}
\begin{align}
    h(x) &=-35.0 \omega^2-29.9 \theta^2 - 5.1 \omega \theta + 1.3 \omega + 7.3 \theta + 12.0, \nonumber \\
    \gamma(r) &= r, ~~ r \in \mathbb{R}_{\geqslant 0}. \nonumber
\end{align}
\section{Conclusions and Future Work}
In this paper, we have introduced a class of robust discrete-time Control Barrier Functions (R-DTCBFs) that can be used to synthesize safe controllers. To construct such an R-DTCBF, we have proposed a counterexample-guided approach. As part of our future work, we aim to improve the method to accelerate the synthesis of R-DTCBFs for higher-dimensional systems, since the computational time of the proposed verification algorithm increases exponentially with the number of states.


\bibliographystyle{IEEEtran}
\bibliography{IEEEabrv}

\end{document}